\documentclass[12pt]{amsart}
\usepackage[mathscr]{eucal}
\usepackage{amssymb,amsmath}
\usepackage{amsfonts}
\usepackage{a4}

\newtheorem{theorem}{Theorem}
\newtheorem{proposition}{Proposition}
\newtheorem{lemma}{Lemma}
\newtheorem{corollary}{Corollary}

\begin{document}
\title{A pairing between super Lie-Rinehart and periodic cyclic homology.}

\author[Tomasz Maszczyk]{Tomasz Maszczyk\dag}
\address{Institute of Mathematics\\
Polish Academy of Sciences\\
Sniadeckich 8\newline 00--956 Warszawa, Poland\\
\newline Institute of Mathematics\\
University of Warsaw\\ Banacha 2\newline 02--097 Warszawa, Poland}
\email{maszczyk@mimuw.edu.pl}

\thanks{\dag The author was partially supported by the KBN grant 1P03A 036 26.}
\thanks{{\em Mathematics Subject Classification:} Primary 16E40, 17B35, 19K56, Secondary 46L87.}

\begin{abstract}
We consider a pairing producing various cyclic Hochschild
cocycles, which led Alain Connes to cyclic cohomology. We are
interested in geometrical meaning and homological properties of
this pairing. We define a non-trivial pairing between the homology
of a Lie-Rinehart (super-)algebra with coefficients in some
partial traces and relative periodic cyclic homology. This pairing
generalizes the index formula for summable Fredholm modules, the
Connes-Kubo formula for the Hall conductivity and the formula
computing the $K^{0}$-group of a smooth noncommutative torus. It
also produces new homological invariants of proper maps
contracting each orbit contained in a closed invariant subset in a
manifold acted on smoothly by a connected Lie group. Finally we
compare it with the characteristic map for the Hopf-cyclic
cohomology.
\end{abstract}

\maketitle

\paragraph{\textbf{1. Introduction}} Let $G$ be a simply-connected
Lie group acting smoothly on a smooth manifold $N$ and $Z$ be a
closed invariant submanifold. Let a smooth map $N\rightarrow M$
contract these orbits. In the dual language of algebras of smooth
functions the situation can be described as follows. We have a Lie
algebra $\mathfrak{g}$ acting by derivations on an algebra $B$,
fixing an ideal $J$. We have also a homomorphism of algebras
$\pi^{*} : A\rightarrow B$ such that
$\mathfrak{g}(\pi^{*}A)\subset J$, or equivalently, we have a
homomorphism of algebras
\begin{align}A\rightarrow B\times_{B/J}(B/J)^{\mathfrak{g}}.\end{align} The last
homomorphism of algebras of smooth functions describes a
continuous map of topological spaces
\begin{align}N\sqcup_{P}P/G\rightarrow M.\end{align}

We generalize this construction to the case of families,
parameterized by commutative super-spaces, of noncommutative
super-spaces, with the total space acted by super-Lie-Rinehart
algebras over the base algebra, as follows. Let $(L,R)$ be a
$\mathbb{Z}/2$-graded Lie-Rinehart algebra over a
$\mathbb{Z}/2$-graded-commutative ring $R$ containing rational
numbers, with a subring of constants $k:=H^{0}(L, R; R)=R^{L}$,
acting (from the left) by super-derivations on a
$\mathbb{Z}/2$-graded associative $R$-algebra $B$. Provided a
homomorphism of $\mathbb{Z}/2$-graded associative $k$-algebras
$A\rightarrow B\times_{B/J}(B/J)^{L}$ is given, we prove the
following theorem.
\begin{theorem} There exists a nontrivial canonical $k$-bilinear pairing
\begin{align}H_{p}(L, R; H^{0}(B, (J^{p})^{*}))\otimes_{k} \textbf{\textsl{S}}(HC_{p+2}(A/k))\rightarrow k.\end{align}
\end{theorem}
Here $H_{p}(L, R; - )$ denotes the $p$-th super-Lie-Rinehart
homology \cite{Fu}, \cite{Rin}, $H^{0}(B, - )$ denotes the $0$-th
Hochschild cohomology, $( - )^{*}=Hom_{k}( - , k)$ and
$\textbf{\textsl{S}}: HC_{p+2}(A/k)\rightarrow HC_{p}(A/k))$ is
the periodicity map of Connes on the relative cyclic homology
\cite{Co}. The above theorem implies, after passing to the inverse
limit with respect to $\textbf{\textsl{S}}$, the existence of a
canonical bilinear pairing with the relative periodic cyclic
homology of $A$.
\begin{corollary} There exists a nontrivial canonical bilinear pairing
\begin{align}H_{p}(L, R; H^{0}(B, (J^{p})^{*}))\otimes_{k} HP_{p}(A/k)\rightarrow k.\end{align}
\end{corollary}
The latter pairing induces (is equivalent to, if $k=R$ and $k$ is
a field) the following $k$-linear map
\begin{align}H_{p}(L, R; H^{0}(B, (J^{p})^{*}))\rightarrow HP^{p}(A/k),\end{align}
which can be regarded as a kind of \textit{characteristic map}. We
will compare it with the Connes-Moscovici characteristic map
\cite{CoMo, CoMo1} and with the cup-product  of the second kind of
Khalkhali-Rangipour \cite{KhaRan} in Hopf-cyclic cohomology.

We will show four classes of examples for which our pairing (or
the characteristic map) is known to be, in general, non-trivial
and its values have important geometric interpretations. First one
is the creation of nontrivial homology classes by contracting
orbits making sense in classical differential geometry, the second
is the index formula for summable Fredholm modules \cite{Co, Co2}
third is the Connes-Kubo formula for the Hall conductivity in the
quantum Hall effect \cite{CHMM1, Bel+E+S, Co, Mc, MM, Xia}, and
the fourth computes $K^{0}$-group of a noncommutative torus in
terms of characteristic numbers of smooth Powers-Rieffel
projections \cite{Co1, PimVoi, Rief}.

Analogous considerations give us the following ``dual" variant,
seemingly more fundamental, of our construction for $((L, R),
J\subset B)$ and $A\rightarrow B$ as above.
\begin{theorem} There exists a nontrivial canonical $k$-linear ``dual characteristic
map"
\begin{align}HP_{p}(A/k)\rightarrow H^{p}(L, R; H_{0}(B, J^{p})).\end{align}
\end{theorem}

\vspace{3mm}
\paragraph{\textbf{2. Construction}} We consider the bilinear
pairing
\begin{align}C_{p}(L, R; H^{0}(B, (J^{p})^{*}))\otimes_{k} \bigotimes\!^{p+1}_{k}A\rightarrow k,\end{align}
$$(\tau\otimes X_{1}\wedge\ldots\wedge X_{p})\otimes
(a_{0}\otimes\ldots\otimes a_{p})\mapsto
\sum_{\sigma\in\Sigma_{p}} (\pm) \tau
(a_{0}X_{\sigma{(1)}}(a_{1})\ldots X_{\sigma{(p)}}(a_{p})),$$
where the sign is determined uniquely by the convention of
transposition of homogeneous symbols from the left hand side to
the position on the right hand side. In the sequel we will use
homogeneous elements, the above sign convention and the
abbreviated notation $X = X_{1}\wedge\ldots\wedge X_{p}\in
\bigwedge^{p}_{R}L$. Let $\tau\in H^{0}(B, (J^{p})^{*})=(J^{p}/[B,
J^{p}])^{*}$. The latter space is a right $(L,R)$-module. The
super-Lie-Rinehart boundary operator $\partial: C_{p}(L, R; -
)\rightarrow C_{p-1}(L, R; - )$, where $C_{p}(L, R; - )=( -
)\otimes_{R}\bigwedge^{p}_{R}$, computing homology with values in
$(L,R)$-modules is an obvious minimal common generalization of the
super-Lie boundary operator from\cite{Fu} and the Lie-Rinehart
boundary operator from \cite{Rin}. By $Z_{p}(L, R; -)$ (resp.
$B_{p}(L, R; -)$) we denote cycles (resp. boundaries) in this
complex. By \textbf{\textsl{b}} (resp. \textbf{\textsl{t}},
\textbf{\textsl{B}}) we denote the Hochschild boundary (resp.
cyclic operator, Connes B-operator) used in cyclic homology
\cite{Co1}. In the lemmas below we apply the above pairing to
various pairs of submodules of super-Lie-Rinehart and Hochschild
chains.
\begin{lemma}
\begin{align}C_{p}(L, R; H^{0}(B, (J^{p})^{*}))\cdot {\rm im} ( \textbf{\textsl{b}} ) = 0.\end{align}
\end{lemma}
\textit{Proof.}
\begin{align}(\tau\otimes X)\cdot \textbf{\textsl{b}} (a_{0}\otimes\ldots\otimes a_{p})=0.\end{align}
$\Box$
\begin{lemma}
\begin{align}Z_{p}(L, R; H^{0}(B, (J^{p})^{*}))\cdot {\rm im}(1-\textbf{\textsl{t}}) = 0.\end{align}
\end{lemma}
\textit{Proof.}
\begin{align}(\tau\otimes X)\cdot (1-\textbf{\textsl{t}})(a_{0}\otimes\ldots\otimes a_{p})= \pm\partial(\tau\otimes X)\cdot (a_{p}a_{0}\otimes a_{1}\otimes\ldots\otimes a_{p-1}).\end{align}
$\Box$

From the last two lemmas we get
\begin{corollary} There exists a canonical bilinear pairing
\begin{align}Z_{p}(L, R; H^{0}(B, (J^{p})^{*}))\otimes HC_{p}(A/k)\rightarrow k.\end{align}
\end{corollary}
One could expect that the above pairing descends to Lie algebra
homology. But it is not true without an appropriate replacement on
the level of cyclic homology.
\begin{lemma}
\begin{align}B_{p}(L, R; H^{0}(B, (J^{p})^{*}))\cdot {\rm ker}(\textbf{\textsl{B}}:
HC_{p}(A/k)\rightarrow HH_{p+1}(A/k))= 0.\end{align}
\end{lemma}
\textit{Proof.} The following formula is an analog of the Stokes
formula
\begin{align}(\tau\otimes X)\cdot \textbf{\textsl{B}}(a_{0}\otimes\ldots\otimes a_{p-1})=\pm p\ \partial(\tau\otimes X)\cdot (a_{0}\otimes\ldots\otimes a_{p-1}).\end{align}
$\Box$

By the long exact sequence of Connes
\begin{align}\ldots\rightarrow HH_{p+2}(A/k)\stackrel{\textbf{\textsl{I}}}{\rightarrow}HC_{p+2}(A/k)\stackrel{\textbf{\textsl{S}}}{\rightarrow}HC_{p}(A/k)\stackrel{\textbf{\textsl{B}}}{\rightarrow}HH_{p+1}(A/k)\rightarrow\ldots \end{align}
we have
$$ker(\textbf{\textsl{B}}: HC_{p}(A/k)\rightarrow HH_{p+1}(A/k))=im (\textbf{\textsl{S}}: HC_{p+2}(A/k)\rightarrow HC_{p}(A/k)).$$
Together with Lemma 3 and Corollary 2 this gives the pairing
\begin{align}H_{p}(L, R; H^{0}(B, (J^{p})^{*}))\otimes_{k} im (\textbf{\textsl{S}}: HC_{p+2}(A/k)\rightarrow HC_{p}(A/k))\rightarrow k\end{align}
desired in Theorem 1.

In order to show that it is non-trivial and interesting we
consider the following classes of examples.

\vspace{3mm}
\paragraph{\textbf{3. Example: Contracting orbits}} Before we
present non-classical examples, we want to explain the classical
case in differential topology. Let $N$ be a compact manifold
(resp. singular with boundary) acted on by a connected Lie group
$G$ with Lie algebra $\mathfrak{g}$, $P$ be a closed invariant
subset (resp. containing the singular locus or boundary) and
$J\subset B:=C^{\infty}(N)$ be a $\mathfrak{g}$-invariant ideal of
smooth functions on $N$ vanishing along $P$. The action of
$\mathfrak{g}$ on differential forms on $N$ is a representation of
a $\mathbb{Z}$-graded super-Lie algebra linearly spanned by
symbols $(d, \iota_{X}, \mathcal{L}_{X})$, where
$X\in\mathfrak{g}$, of degrees $(1, -1, 0)$, subject to the
relations
\begin{align}[d, d] = 0,\  [\iota_{X}, \iota_{Y}]=0,\ [d, \mathcal{L}_{X}]=0,
\end{align}
$$[d, \iota_{X}]=\mathcal{L}_{X},\ [\mathcal{L}_{X},
\iota_{Y}]=\iota_{[X, Y]},\ [\mathcal{L}_{X},
\mathcal{L}_{Y}]=\mathcal{L}_{[X, Y]}.$$ We will use the following
 consequence of these relations

\begin{align}[d, \iota_{X_{1}}\ldots\iota_{X_{p}}]=\end{align}
$$=\sum_{i}(-1)^{i-1}\iota_{X_{1}}\ldots\widehat{\iota_{X_{i}}}\ldots\iota_{X_{p}}\mathcal{L}_{X_{i}} +
\sum_{i<j}(-1)^{i+j-1}\iota_{[X_{i}, X_{i}
]}\iota_{X_{1}}\ldots\widehat{\iota_{X_{i}}}\ldots\widehat{\iota_{X_{j}}}\ldots\iota_{X_{p}}.$$
Every smooth measure $\mu $ on $N\setminus P$, (i.e. a
differential top degree form with values in the orientation
bundle), such that for every element $f\in J^{p}$ the product
$f\mu$ extends to a smooth measure on the whole $N$, defines an
element
\begin{align}\int_{Y}(-)\mu\in H^{0}(B, (J^{p})^{*})=(J^{p})^{*}.\end{align} The right
$\mathfrak{g}$-action on such element  reads as
\begin{align}(\int_{Y}(-)\mu)\cdot X = -\int_{Y}(-)\mathcal{L}_{X}\mu\end{align}

\begin{proposition} If a chain
\begin{align}\sum\int_{Y}(-)\mu\otimes X_{1}\wedge\ldots\wedge X_{p} \in C_{p}(\mathfrak{g}, H^{0}(B, (J^{p})^{*})),\end{align}
where $\mu $'s are smooth measures on $Y$ as above, is a cycle
(resp. a boundary) then the differential form $\sum
\iota_{X_{1}}\ldots \iota_{X_{p}}\mu $ is closed (resp. exact).
\end{proposition}

\textit{Proof.} The cycle condition for our chain
\begin{align}\partial\sum\int_{Y}(-)\mu\otimes X_{1}\wedge\ldots\wedge X_{p}=0\end{align}
is equivalent to
\begin{align}\sum\sum_{i}(-1)^{i}X_{1}\wedge\ldots \widehat{X_{i}}\ldots\wedge X_{p}\otimes \mathcal{L}_{X_{i}}\mu\  +\end{align}
$$+
\sum\sum_{i<j}(-1)^{i+j-1}[X_{i}, X_{i} ]\wedge X_{1}\wedge\ldots
\widehat{X_{i}}\ldots \widehat{X_{j}}\ldots\wedge X_{p}\otimes
\mu=0,$$ which implies that
\begin{align}\sum\sum_{i}(-1)^{i}\iota_{X_{1}}\ldots \widehat{\iota_{X_{i}}}\ldots\iota_{X_{p}}\mathcal{L}_{X_{i}}\mu\  +\end{align}
$$+
\sum\sum_{i<j}(-1)^{i+j-1}\iota_{[X_{i}, X_{i} ]}
\iota_{X_{1}}\ldots \widehat{\iota_{X_{i}}}\ldots
\widehat{\iota_{X_{j}}}\ldots \iota_{X_{p}}\mu=0,$$ which is
equivalent to
\begin{align}\sum[d, \iota_{X_{1}}\ldots \iota_{X_{p}}]\mu =0.\end{align}
Since $\mu$ is a top degree form $d\mu=0$, which gives finally
\begin{align}d\sum \iota_{X_{1}}\ldots \iota_{X_{p}}\mu =0. \end{align}
The proof of the implication ``\textit{boundary} $\Rightarrow$
\textit{exact}" is similar. $\Box$

Let us consider now a smooth map $\pi: N\rightarrow M$ into a
compact manifold (resp. singular variety, with boundary) $M$,
contracting each orbit contained in $P$ to a point. Let
$A:=C^{\infty}(M)$ be an algebra of smooth functions on $M$. Then
$\mathfrak{g}(\pi^{*}A)\subset J\subset B$.

Comparing with the canonical map from De Rham homology of currents
to periodic cyclic cohomology we see that our characteristic map
associates with the Lie homology class of the above cycle  a
homology class of the closed current $j$, where
\begin{align}j(\omega):= \sum\int_{N}(\pi^{*}\omega)(\tilde{X}_{1},\ldots , \tilde{X}_{p})\mu = \pm\int_{N}(\pi^{*}\omega)\wedge\sum \iota_{X_{1}}\ldots \iota_{X_{p}}\mu
.\end{align} Here by $\tilde{X}$ we mean the vector field
corresponding to an element $X\in \mathfrak{g}$.

Note that if $\sum \iota_{X_{1}}\ldots \iota_{X_{p}}\mu$ extends
to the whole $N$ then we get the push-forward of the homology
class of the closed current
\begin{align} \left[ \sum\int_{N}(-)\mu\otimes X_{1}\wedge\ldots\wedge X_{p}\right] \mapsto\pm \pi_{!}\left[ \sum\int_{N}(-)\iota_{X_{1}}\ldots \iota_{X_{p}}\mu\right]
\end{align} in periodic cyclic cohomology.

Take for example $N=M=G=S^{1}$ and $\pi={\rm id}:S^{1}\rightarrow
S^{1}$, $P=\emptyset$, $\mu$ the Haar measure and $X\in
\mathfrak{g}$ normalized so that $\iota_{X}\mu = 1$. Then
$\int_{S^{1}}(-)\mu\otimes X $ is a cycle and the characteristic
map gives
\begin{align}\left[ \int_{S^{1}}(-)\mu\otimes X\right] \mapsto \left[ \int_{S^{1}}\right] ,\end{align}
i.e. the fundamental class of $S^{1}$. Though it is a nontrivial
example, it is not very enlightening. Therefore we need more
complicated example to show the point. This time $\sum
\iota_{X_{1}}\ldots \iota_{X_{p}}\mu$ will not extend to the whole
$N$, so the push-forward of the respective current will not be
defined. However, our characteristic map still will define a
nontrivial homology class on $M$. To see this, let us take a
cylinder $N=S^{1}\times [-\pi, \pi]$ with coordinates $(\varphi,
\psi)$ ($\varphi$ - circular coordinate, $\psi\in [-\pi, \pi]$).
Consider the following smooth action of the additive Lie group
$G=\mathbb{R}$ on $N$
\begin{align*}
t\cdot (\varphi, \psi):=\left\{ \begin{array}{cc}
                          (\varphi+t,2\arctan(\tan\frac{\psi}{2}
                          +t)) & {\rm if} \psi\neq \pm \pi,\\
                          (\varphi+t,\pm \pi)) & {\rm if} \psi= \pm
                          \pi.
                                                    \end{array}\right.
                                                    \end{align*}
It is generated by a vector field
\begin{align*}
X=\frac{\partial}{\partial \varphi}+(1+\cos \psi
)\frac{\partial}{\partial \psi}.
\end{align*}
Let us take $P=\partial N$, which is the union of compact orbits
of the above action and the ideal $J =(1+\cos \psi )$ in the
algebra $B=C^{\infty}(N)$, vanishing along $P$. The form
\begin{align*}
\mu=\frac{1}{2\pi}\frac{d\psi\wedge d\varphi}{1+\cos \psi}
\end{align*}
is defined on $N\setminus P$ and invariant. We have
\begin{align*}
\iota_{X}\mu=\frac{1}{2\pi}d\varphi
-\frac{1}{2\pi}\frac{d\psi}{1+\cos \psi},
\end{align*} \vspace{3mm}
which does not extend onto the whole $N$.

Take now a subvariety $M\subset \mathbb{R}^{3}$
\begin{align*}
z^{2}
=\varepsilon^{2}(\frac{x}{\sqrt{x^{2}+y^{2}}}+1)^{2}-(\sqrt{x^{2}+y^{2}}-1)^{2},\
\ \ (x, y)\neq (0, 0),
\end{align*}
where a parameter $\varepsilon \in (0, 1/2)$, homeomorphic to a
torus with one basic  cycle contracted to the unique singular
point $(-1, 0, 0)$, and a smooth map $\pi:N\rightarrow M$ of the
form
\begin{align*}
\pi^{*}x & =\cos \psi(1+\varepsilon(1+\cos \psi)\cos \varphi),\\
\pi^{*}y & =\sin \psi(1+\varepsilon(1+\cos \psi)\cos \varphi),\\
\pi^{*}z & =\varepsilon(1+\cos \psi)\sin\varphi.
\end{align*}
One can check that $X(\pi^{*}x), X(\pi^{*}y), X(\pi^{*}z)\in J$
which implies that for $A=C^{\infty}(M)$
$$X(\pi^{*}A)\subset J.$$
We have $H_{1}(N, \mathbb{Z})=\mathbb{Z}$ generated by the
homology class of one boundary circle $(\psi=\pi)$, $H_{1}(M,
\mathbb{Z})=\mathbb{Z}$ generated by the homology class of the
ellipse $(x^{2}+y^{2}=1, z=\varepsilon(x+1))$. The topology of the
map $\pi:N\rightarrow M$ is following. It contracts the boundary
circles of the cylinder $N$ to the singular point of $M$.
Therefore it kills the generator of $H_{1}(N, \mathbb{Z})$. But it
also creates the generator of $H_{1}(M, \mathbb{Z})$. The killing
property of $\pi$ is described by the nullity of the induced map
$H_{1}(N, \mathbb{Z})\rightarrow H_{1}(M, \mathbb{Z})$. We will
show that the creating property of $\pi$ is described by our
characteristic map. Let
$$\omega=\frac{1}{2\pi}\frac{xdy-ydx}{x^{2}+y^{2}}$$
be a closed 1-form on $M$, whose period over the generator of
$H_{1}(M, \mathbb{Z})$ is equal to 1. We have
$\pi^{*}\omega=d\psi$. Let us compute our pairing of the Lie
homology class with the De Rham cohomology class of $\omega$
$$\left[ \int_{N}(-)\mu\otimes X\right]\cdot[\omega]=-\int_{N}(\pi^{*}\omega)\wedge\iota_{X}\mu=-\frac{1}{4\pi^{2}}\int_{N}d\psi\wedge d\varphi=1.$$
Therefore our characteristic map applied to $\left[
\int_{N}(-)\mu\otimes X\right]$ gives the homology class of the
current homological to the period over the generator of $H_{1}(M,
\mathbb{Z})$. \vspace{3mm}

\paragraph{\textbf{4. Example: Index formula}} Let us assume that $k=\mathbb{C}$ and we have a $p$-summable
even involutive Fredholm module $(A, H, F)$ \cite{Co, Co2}, i.e.
$A$ is a $\mathbb{Z}/2$-graded $*$-algebra, $H$ is a
$\mathbb{Z}/2$-graded Hilbert space with a grading preserving
$*$-representation $A\rightarrow B(H)$, and $F$ is an odd
self-adjoint involution on $H$ such that
\begin{align}[F, A]\subset L^{p}(H),\end{align}
where $L^{p}(H)$ denotes the $p$-th Schatten ideal in $B(H)$.

Let us define now a $\mathbb{Z}/2$-graded abelian super-Lie
algebra generated by one odd element $d$
\begin{align}\mathfrak{g}:=\mathbb{C}\cdot d,\end{align}
a $\mathbb{Z}/2$-graded associative algebra
\begin{align}B:=B(H)\end{align}
and an ideal
\begin{align}J:=L^{p}(H).\end{align}
The formula
\begin{align}db:=[F, b]\end{align}
defines the left action of $\mathfrak{g}$ on $B$ by derivations
and obviously $J$ is a $\mathfrak{g}$-ideal.

The projection into the first cartesian factor defines an
isomorphism of $\mathbb{Z}/2$-graded associative algebras
\begin{align} B\times_{B/J}(B/J)^{\mathfrak{g}}\stackrel{\cong}{\rightarrow}\{ b\in B(H) \mid [F, b]\in L^{p}(H)\}\end{align}
and a $*$-homomorphism of $\mathbb{Z}/2$-graded associative
$C^{*}$-algebras $A\rightarrow B\times_{B/J}(B/J)^{\mathfrak{g}}$
is equivalent to a structure of $p$-summable even involutive
Fredholm module $(A, H, F)$. By functoriality of Connes' long
exact sequence it is enough to consider our pairing for $A:=
B\times_{B/J}(B/J)^{\mathfrak{g}}$.

Since the super-trace is a linear functional on the $p$-th power
of the ideal $J:=L^{p}(H)$ which vanishes on super-commutators and
the super-Lie algebra $\mathfrak{g}=\mathbb{C}\cdot d$ is abelian,
the element ${\rm str} \otimes d\wedge\ldots\wedge d\in H^{0}(A,
(J^{p})^{*})\otimes \bigwedge^{p}\mathfrak{g}$ is a cycle. On the
other hand, for any even self-adjoint idempotent $e\in A$ the
element $e\otimes\ldots \otimes e\in \bigotimes^{p+1}A$ is a
cyclic cycle for even $p$, which is in the image of the
periodicity operator $\textbf{\textsl{S}}$. We can compute our
pairing of homology classes of these two cycles which gives the
index of a Fredholm operator
\begin{align}[ {\rm str}\otimes d\wedge\ldots\wedge d ]\cdot [e\otimes\ldots \otimes e] = c_{p}\ {\rm Index} (e_{11}F_{01}e_{00}),\end{align}
in general a non-zero number. Here $F_{01}: H_{0}\rightarrow
H_{1}$ (resp. $e_{00}: H_{0}\rightarrow H_{0}$, $e_{11}:
H_{1}\rightarrow H_{1}$) is a unitary block of $F$ (resp.
self-adjoint idempotent block of $e$) under the orthogonal
decomposition $H=H_{0}\oplus H_{1}$ into even and odd part.

\vspace{3mm}
\paragraph{\textbf{5. Example: Connes-Kubo formula}} Let $\mathfrak{g}$ be an
abelian Lie algebra and $A=B=J$. If $\tau$ is a
$\mathfrak{g}$-invariant trace on $A$ then this is obvious that
for all $X_{1},\ldots, X_{p}\in \mathfrak{g}$ the chain
\begin{align}\tau\otimes X_{1}\wedge\ldots\wedge X_{p} \end{align}
is a cycle hence defines a homology class. This construction is
next adapted to the geometry of the Brillouin zone. Its pairing
with an appropriately normalized even dimensional class $[e\otimes
e \otimes e]$ in $HP_{2}(A)$ computes the Hall conductivity
$\sigma$ in noncommutative geometric models of quantum Hall effect
\begin{align}[\tau\otimes \sum_{i=1}^{g}X_{i}\wedge X_{i+g}]\cdot [e\otimes e\otimes e] = \sigma ,\end{align}
in general a non-zero integer \cite{CHMM1, Bel+E+S, Co, Mc, Xia}
or rational number \cite{MM}, depending on the model.

\vspace{3mm}
\paragraph{\textbf{6. Example: $K_{0}$ of a noncommutative torus}} Formally it is the same construction
as in the previous example adapted to the context of
non-commutative geometry of the noncommutative 2-torus
\cite{Rief}. Let $A$ be the dense subalgebra of ``smooth functions
on the noncommutative 2-torus" \cite{Co1} of the $C^{*}$-algebra
generated by two unitaries $U,V$ subject to the relation
\begin{align}UV=e^{2\pi i\theta}VU\end{align} with an irrational
real $\theta$. The Lie group $S^{1}\times S^{1}$ acts on $A$ by
automorphisms and its Lie algebra $\mathfrak{g}$ spanned by
commuting elements $X, Y$ acts by derivations such that
\begin{align}X(U)=2\pi i U,\ X(V)=0,\end{align}
\begin{align}Y(U)=0,\ Y(V)=2\pi i V.\end{align}
Every element $a\in A$ can be uniquely expanded as $a=\sum
a_{mn}U^{m}V^{n}$. Then the functional $\tau (a)=a_{00}$ is a
$\mathfrak{g}$-invariant trace. Again, we have a homology class
$[\tau\otimes X\wedge Y]$. It is known that $K_{0}(A)=\mathbb{Z}
+\mathbb{Z}\cdot\theta\subset\mathbb{C}$ where the identification
is done by this trace \cite{PimVoi, Rief}. Any selfadjoint
idempotent $e\in A$ is determined by its trace $\tau
(e)=p-q\cdot\theta$ uniquely up to unitary equivalence. This is
our pairing in dimension zero
\begin{align}[\tau]\cdot [e]=p-q\cdot\theta.\end{align}
Our pairing in dimension two computes the number $q$
\begin{align}[\tau\otimes X\wedge Y]\cdot [e\otimes e\otimes e] = q\cdot 2\pi i .\end{align}
This means that our pairings, defined \textit{a priori} over
$\mathbb{C}$, detect fully the $K_{0}$-group isomorphic to
$\mathbb{Z}\oplus \mathbb{Z}$.

\vspace{3mm}
\paragraph{\textbf{7. Comparison with other constructions.}} In \cite{CoMo,
CoMo1} the following pairing (Connes-Moscovici characteristic map)
is considered
\begin{align}HP^{p}_{(\delta , \sigma)}(H)\otimes Tr_{(\delta , \sigma )}(A)\rightarrow HP^{p}(A)\end{align}
for any Hopf algebra $H$ with $(\delta , \sigma)$ a modular pair
in involution, acting on an algebra $A$, where by $Tr_{(\delta ,
\sigma )}(A)$ we denote the space of $(\delta , \sigma )$-traces
on $A$. Taking $H=U(\mathfrak{g})$, $\delta =\epsilon $, $\sigma
=1$ one has \cite{CoMo1}
\begin{align}HP^{p}_{(\epsilon , 1)}(U(\mathfrak{g} ))=\bigoplus_{i\equiv p\ (mod\ 2)
}H_{i}(\mathfrak{g} ),\end{align}
\begin{align}Tr_{(\epsilon , 1 )}(A)=H^{0}(A, A^{*})^{\mathfrak{g} }.\end{align} Then we have the
following commuting diagram
\begin{align}\begin{array}{ccc}
H_{p}(\mathfrak{g})\otimes H^{0}(A,A^{*})^{\mathfrak{g}} &
\rightarrow & H_{p}(\mathfrak{g}, H^{0}(A, A^{*})) \\
\downarrow &  & \downarrow \\
HP^{p}_{(\epsilon , 1)}(U(\mathfrak{g} ))\otimes Tr_{(\epsilon , 1
)}(A) & \rightarrow & HP^{p}(A),
\end{array}\end{align}
where left vertical and upper horizontal arrows are canonical, the
bottom horizontal arrow is the Connes-Moscovici characteristic map
and the right vertical arrow is our characteristic map for
$A=B=J$. The main difference between these two characteristic maps
is the position of traces: in the Connes-Moscovici map traces are
paired with cyclic periodic cohomology while in our map they are
coefficients of Lie algebra homology.

Recently \cite{KhaRan} a new pairing with values in cyclic
cohomology (the cup product of the second kind)
\begin{align}HC^{p}_{H}(C, M)\otimes HC^{q}_{H}(A, M)\rightarrow HC^{p+q}(A)\end{align}
has been presented, which allows to consider in this pairing
cyclic cohomology with nontrivial coefficients in the sense of
\cite{HaKhaRaSo}. It is defined for a Hopf algebra $H$, an
$H$-module algebra $A$, an $H$-comodule algebra $B$, an $H$-module
coalgebra $C$ acting on $A$ in a suitable sense and any stable
anti-Yetter-Drinfeld (SAYD) module $M$ over $H$. For
$C=H=U(\mathfrak{g} )$, $q=0$ and $M=k_{(\epsilon , 1)}$ trivial
one dimensional SAYD-module one gets again the Connes-Moscovici
characteristic map.

Since $U(\mathfrak{g})$ is a cocommutative Hopf algebra, any
$U(\mathfrak{g} )$-module $M$ with a trivial $U(\mathfrak{g}
)$-comodule structure is a SAYD-module. Taking $C=H=U(\mathfrak{g}
)$, $q=0$ and $M=H^{0}(A,A^{*})$ one has the Khalkhali-Rangipour
cup product of the second kind
\begin{align}HC^{p}_{U(\mathfrak{g} )}(U(\mathfrak{g}), H^{0}(A,A^{*}))
\otimes HC^{0}_{U(\mathfrak{g} )}(A, H^{0}(A,A^{*}))\rightarrow
HC^{p}(A).\end{align}

The trace evaluation map $H^{0}(A,A^{*})\otimes A\rightarrow k$
defines a distinguished element in $HC^{0}_{U(\mathfrak{g} )}(A,
H^{0}(A,A^{*}))$ and consequently the following characteristic map
\begin{align}HC^{p}_{U(\mathfrak{g} )}(U(\mathfrak{g}), H^{0}(A,A^{*}))
\rightarrow HC^{p}(A),\end{align} by taking the above cup product
with this distinguished element. In fact this map comes from the
morphism of cyclic objects, so it can be pushed to the periodic
cyclic cohomology, hence we get a map \begin{align}
HP^{p}_{U(\mathfrak{g} )}(U(\mathfrak{g}), H^{0}(A,A^{*}))
\rightarrow HP^{p}(A). \end{align} By Theorem 5.2 of \cite{JarSte}
(note that in \cite{JarSte} authors use cyclic objects related to
cyclic objects from \cite{HaKhaRaSo} by Connes's cyclic duality
\cite{KhaRan'}, transforming homology into cohomology) one has
\begin{align}HP^{p}_{U(\mathfrak{g} )}(U(\mathfrak{g}), H^{0}(A,A^{*}))=\bigoplus_{i\equiv p\ (mod\ 2)
}H_{i}(\mathfrak{g}, H^{0}(A,A^{*})).\end{align} In this case, our
characteristic map factorizes through (51)
\begin{align}\begin{array}{ccc}
            &                    & H_{p}(\mathfrak{g}, H^{0}(A, A^{*})) \\
            &     \swarrow       & \downarrow \\
HP^{p}_{U(\mathfrak{g} )}(U(\mathfrak{g}), H^{0}(A,A^{*}))
 & \rightarrow & HP^{p}(A),
\end{array}\end{align}
where the south-west arrow is an embedding onto a direct summand
in the decomposition (52).

However, in general, there is no way to extend a partial trace
from the ideal $J^{p}$ to the trace defined on the whole algebra,
hence there is no a canonical element to pair with as in (50).
Therefore the characteristic map \textit{\'{a} la}
Khalkhali-Rangipour $``HP^{p}_{U(\mathfrak{g} )}(U(\mathfrak{g}),
H^{0}(B,(J^{p})^{*}))
  \rightarrow  HP^{p}(A)"$ is not defined in general. In particular, the index pairing discussed in the
  paragraph 4 cannot be obtained in this way.

In spite of this discrepancy we expect that, after appropriate
modifications of Hopf-cyclic cohomology (or its extended version
\cite{KhaRan''} working in the case of enveloping algebras of
Lie-Rinehart algebras) with appropriate coefficients, our
construction could be generalized to (super, extended) Hopf-cyclic
cohomology. The crucial property this generalization should
satisfy is the above index pairing.


\begin{thebibliography}{9999999}

\bibitem{Bel+E+S} Bellissard, J.; van Elst, A.; Schulz-Baldes, H.: The
non-commutative geometry of the quantum Hall effect. {\em J. Math.
Phys.} {\bf 35} (1994), 5373-5451.

\bibitem{Mc}  McCann, P. J.: Geometry and the integer quantum Hall effect,
in {\em Geometric Analysis and Lie Theory in Mathematics and
Physics.} pp 132-208, Edited by A.L. Carey and M.K. Murray
Cambridge Univ. Press, Cambridge 1998.

\bibitem{CHMM1}  Carey, A.;  Hannabus, K.;  Mathai, V.;   McCann, P.:
Quantum Hall Effect on the hyperbolic plane. {\em Commun. Math.
Physics.} {\bf 190}, No. 3  (1976), 629-673.

\bibitem{Co1} Connes, A.: $C^{*}$-alg\'{e}bres et g\'{e}om\'{e}trie diff\'{e}rentielle.
{\em C. R. Acad. Sci. Paris} Ser. A-B,{\bf 290}, 1980.

\bibitem{Co} Connes, A.: Noncommutative differential geometry.
{\em Publ. Math. I.H.E.S.} {\bf 62} (1986), 257-360.

\bibitem{Co2} Connes, A.: Noncommutative geometry. Acad. Press, Inc., San
Diego, CA, (1994).

\bibitem{CoMo} Connes, A.; Moscovici, H.: Cyclic cohomology and Hopf algebras.
{\em Lett. Math. Phys.} {\bf 48} (1999), 97-108.

\bibitem{CoMo1} Connes, A.; Moscovici, H.: Hopf algebras, cyclic cohomology and the transverse index
theorem. {\em Comm. Math. Phys.} {\bf 198} (1998), 199-246.

\bibitem{HaKhaRaSo} Hajac, P. M.; Khalkhali, M.; Rangipour, B.;  and Y.
Sommerh\"{a}user, Y.: Hopf-cyclic homology and cohomology with
coefficients. {\em C. R. Math. Acad. Sci. Paris} \textbf{338}
(2004), No. 9, 667-672.

\bibitem{JarSte} Jara, P.;  Stefan, D.: Cyclic homology of Hopf
Galois extensions and Hopf algebras. {\em preprint, arXiv:
math.KT/0307099}.

\bibitem{KhaRan''} Khalkhali, M.;  Rangipour, B.: Cyclic cohomology of (extended)
Hopf algebras. Noncommutative geometry and quantum groups (Warsaw,
2001), 59-89, {\em Banach Center Publ.}, \textbf{61}, 2003.

\bibitem{KhaRan'} Khalkhali, M.;  Rangipour, B.: A note on cyclic duality and Hopf algebras. {\em Comm. Alg.} {\bf
33}, No. 3, (2005), 763-773.

\bibitem{KhaRan} Khalkhali, M.;  Rangipour, B.: Cup Products in Hopf-Cyclic
Cohomology. {\em C. R. Acad. Sci. Paris, Ser.} I. {\bf 340}
(2005), 9-14.

\bibitem{Fu} Leites, D.A.; Fuks, D.B.: Cohomology of  Lie
superalgebras, {\em Dokl. Bolg. Akad. Nauk}, {\bf 37}, No. 10
(1984), 1294-1296.

\bibitem{MM}  Marcolli, M.; Mathai, V.: Twisted Higher Index Theory on Good
Orbifolds, II:  Fractional Quantum Numbers. {\em  Comm. Math.
Phys.}, {\bf 201} (2001), No. 1, 55-87.


\bibitem{PimVoi}  Pimsner, M.; Voiculescu, D.: Exact sequences for $K$ groups
and $Ext$ groups of certain cross products $C^{*}$-algebras, {\em
J. Operator Theory}, {\bf 4} (1980), 93-118.

\bibitem{Rief}  Rieffel, M.: $C^{*}$-algebras associated with irrational rotations, {\em Pac.
J. Math.} {\bf 93} (1981), 415-429.

\bibitem{Rin}  Rinehart, G.: Differential forms on general commutative algebras,
{\em Trans. Amer. Math. Soc.} {\bf 108} (1963), 195-222.


\bibitem{Xia}  Xia, J.: Geometric invariants of the quantum Hall effect,
Commun. {\em Math. Phys.} {\bf 119} (1988), 29-50.



\end{thebibliography}
\end{document}